\documentclass[12pt]{amsart}
\usepackage{amsmath}
\usepackage{amsthm}
\usepackage{amsfonts}
\usepackage{graphicx}
\usepackage{amssymb} 
\newenvironment{nouppercase}{
  
  \renewcommand{\uppercasenonmath}[1]{}}{}

\begin{document}

\title[Exact algebraic M(em)brane solutions]
{Exact algebraic M(em)brane solutions}
\author[Jens Hoppe]{Jens Hoppe}
\address{Braunschweig University, Germany}
\email{jens.r.hoppe@gmail.com}

\begin{abstract}
Three classes of new, algebraic, zero-mean-curvature hypersurfaces in pseudo-Euclidean spaces are given.
\end{abstract}

\begin{nouppercase}
\maketitle
\end{nouppercase}
\thispagestyle{empty}
\noindent 
While everyone thinking about relativistic membrane dynamics will easily find a collapsing round sphere of radius $r(t)$ (= some elliptic function), and perhaps (cp. \cite{JH82}) 
\begin{equation}\label{eq1} 
t-z   = \left( \frac{x^2 + y^2}{R^2}-1 \right) \frac{1}{2}R\dot{R} + h\left( \frac{t+z}{2} \right),
\end{equation}
$\dot{h} = \frac{\dot{R}^2}{2}$,\quad $\frac{1}{2}\dot{R}^2 + \frac{1}{4}R^4 = $ const.,  or \cite{JH95} 
\begin{equation}\label{eq2}
\begin{split} 
\mathcal{P}(x)\mathcal{P}(y)\mathcal{P}(t) & = \mathcal{P}(z)\\
\mathcal{P}'^2 & = 4\mathcal{P} (\mathcal{P}^2-1), 
\end{split}
\end{equation}
not a single rotationally non-invariant exact time-like solution in 4-dimensional Minkowski space, arising from the motion of a smooth closed 2-dimensional surface, has been described during the 60 years since Dirac attempted to explain the muon as an excited state of the electron \cite{D62}.
In this note I would like to present 3 classes of new, algebraic, zero-mean-curvature hypersurfaces in ($\mathbb{R}^D, \eta = (\eta_{\mu\nu}) = \text{diag}(1, \pm 1, \ldots, \pm 1, -1)$).
Consider
\begin{equation}\label{eq3}
\begin{split} 
\chi(x) & := (\alpha \cdot x)^n (x \circ x) =: \psi \cdot \phi \\
x \circ x & := x^{\mu}g_{\mu\nu}x^{\nu}, \; (g_{\mu\nu}) = \text{diag}(1, \varepsilon_a, -1) , |\varepsilon_a| = 1\\
\alpha_{\mu} & = (1, 0, \ldots, 0,1)\\
a& = 1,\ldots,D-2 =:M\\
\mu, \nu & = 0,1,\ldots , D-1.
\end{split}
\end{equation}
Straightforwardly one finds
\begin{equation}\label{eq4}
\begin{split} 
\square \chi & := \eta^{\mu\nu}\partial_{\mu}\partial_{\nu}\chi = (4n + 2\gamma)\psi \quad \gamma:=\eta^{\mu\nu}g_{\mu\nu}\\
(\partial \chi)^2 & := \partial_{\mu}\chi\eta^{\mu\nu}\partial_{\nu}\chi = 4\psi^2x^2 + 4n\psi^2\phi.
\end{split}
\end{equation}
The mean curvature of the hypersurfaces $\sum_{D-1} := \lbrace x | \chi(x) = C\rbrace$ vanishes if
\begin{equation}\label{eq5}
\begin{split}
\frac{1}{2} \partial^{\mu}&\chi\partial_{\mu}((\partial \chi)^2)  - (\partial \chi)^2 \square \chi \\
 & = (4n(4n+2\gamma)-(12n^2 + 4n+8))\psi^3\phi +8\gamma\psi^3x^2 
\end{split}
\end{equation}
is zero. If $g_{\mu\nu} \neq \eta_{\mu\nu}$ $\gamma$ must vanish, and $4(n^2-n-2) = 0$, i.e. $n= + 2$ or $-1$. If $g_{\mu\nu}= \eta_{\mu\nu}$ the 2 terms can combine (and $\gamma = D$), resulting in $4(n^2 +(2\gamma-1)n+2(\gamma-1))=0$, hence $n=-1$ or $n = -2(D-1)$. 
Forgetting about $n=-1$ (linear shifts of the coordinates are always allowed anyway) the two classes of solutions thus obtained are
\begin{equation}\label{eq6}
\begin{split}
x^{\mu}x_{\mu} & = (t+z)^{2(D-1)}\cdot C \\
(= x^{\mu}\eta_{\mu \nu}x^{\nu} & = t^2 + \sum^{D-2}_{a=1}\eta_a(x^a)^2 - z^2)
\end{split}
\end{equation}
and 
\begin{equation}\label{eq7}
\begin{split}
\phi(x) & := x^{\mu}g_{\mu\nu}x^{\nu} = \frac{C}{(t+z)^2}\\ 
& (= t^2+ \sum^{D-2}_1 \varepsilon_a(x^a) - z^2)
\end{split}
\end{equation}
where the $\varepsilon_a$, compared with the $\eta_a$ $(a=1 \ldots M=D-2)$ must have 2 more signs differing compared to those agreeing, (i.e. in this case $D$ must be even). Let me discuss the following (simplest, non-trivial) examples, in ordinary Minkowski space $\mathbb{R}^{1,D-1}$ in some detail:
\begin{equation}\label{eq8}
\begin{split}
u(x^{\mu}) & := (t^2 - r^2 -z^2) - C(t+z)^{2M+2} \equiv 0 \\
r^2 &:= x_1^2+ \ldots + x_{D-2 = M}^2 
\end{split}
\end{equation}  
and 
\begin{equation}\label{eq9}
v(x^{\mu}) := \psi \cdot \phi -C' = (t+z)^2\cdot(t^2+r^2-s^2-z^2)- C' \equiv 0
\end{equation}
where $r^2 = x_1^2 + \ldots + x^2_{q=\frac{D-2}{2}+1}$, $s^2 = x^2_{q+1}+\ldots+ x^2_{D-2}$
(and later only the simplest case, $D=4$, $s^2=0$).
Apart from the light-like line $r = 0 = t+z$ (containing the only singular point, $x^{\mu} =0$) (\ref{eq8}) is space-like if $C < 0$, and time-like if $C>0$ (to which we now restrict; for (\ref{eq9}), which for $C'\neq 0$ is regular, we will later take $C' < 0$, corresponding clearly to time-like; the case $C'>0$ is somewhat more subtle). Writing (\ref{eq8}) as 
\begin{equation}\label{eq10}
\begin{split}
r^2 & = \kappa(2t-f(\kappa)) \\
z & = -t+\kappa\\
f & = f_M(\kappa) := \kappa+C\kappa^{2M+1}
\end{split}
\end{equation}   
suitably parametrizes $t$-dependent $M$-dimensional `axially symmetric' surfaces that are compact and convex (due to $f$ being strictly increasing; $\kappa \in [0, \kappa_M(t)]$),
$f(\kappa_M) = 2t > 0$, or, for $t<0$: $\kappa \in [-\kappa_M, 0]$; note the $t \rightarrow -t$, $z\rightarrow -z$ hence $\kappa\rightarrow -\kappa$ invariance of (\ref{eq8}), (\ref{eq9}) resp. (\ref{eq6}), (\ref{eq7})).  
As in the literature on relativistic extended objects, the dynamics is usually given assuming an orthonormal parametrization/gauge (ONG) it would be interesting to reparametrize (\ref{eq10}) by $\kappa = \kappa(t, \varphi)$ such that the hypersurface described by $r(t, \varphi)$ and $z(t,\varphi)$ moves orthogonal to itself, i.e. satisfying $\dot{r}r' + \dot{z}z' = 0$.
Differentiating (\ref{eq10}) with respect to $t$ and $\varphi$ one can then derive
\begin{equation}\label{eq11}
\dot{\kappa} = \frac{2\kappa(2t-f+\kappa f')}{(2t-f-\kappa f')^2 + 4\kappa(2t-f)} =: g(t,\kappa),
\end{equation}   
which on the one hand is reassuring (as $g(t,\kappa) \geqslant 0$ on $\mathbb{R}_+ \times [0, \kappa_M]$ resp. $\mathbb{R}_- \times [-\kappa_M,0]$ and $= 0$ only for $\kappa = z+t =0$), but also showing that even in the string case where everything is always assumed to be of integrable nature, the ODE (\ref{eq11}) does not appear to be easily solvable in explicit terms (let alone the other part of ONG, $\dot{r}^2+\dot{z}^2 + r^M\frac{(r'^2+z'^2)}{\rho^2} = 1$, yielding --using (\ref{eq11})-- an ODE of the form $\frac{\kappa'}{\rho}= h(t,\kappa)$).\\[0.15cm]
Interestingly, the assumption $|\varepsilon_a| = 1$ in (\ref{eq3}) was/is actually not necessary, i.e. not the only one that works. Supposing
\begin{equation}\label{eq12}
\lambda_{\mu \nu} := g_{\mu \mu'}\eta^{\mu' \nu'}g_{\nu'\nu} = \lambda\eta_{\mu \nu} + (1-\lambda)g_{\mu \nu}
\end{equation}  
$\square \chi$ will still be $=(4n + 2\gamma)\psi$, but 
\begin{equation}\label{eq13}
(\partial_{\mu} \chi)^2 = 4\big( (n+1-\lambda)\psi^2\phi + \lambda \psi^2x^2 \big),
\end{equation}  
hence
\begin{equation}\label{eq14}
\begin{split}
\frac{1}{2}\partial^{\mu}\chi \partial_{\mu}\big( (\partial \chi)^2 \big) & = \big( \phi \partial^{\mu}\psi + 2 \psi \eta^{\mu\nu}g_{\nu \rho}x^{\rho}\big) \cdot 4 \big(\lambda\psi\partial_{\mu}\psi x^2 \\
 & \quad\; + \lambda\psi^2 \eta_{\mu\nu}x^{\nu} +(n+1-\lambda)\psi (\phi\partial_{\mu}\psi +g_{\mu\nu}x^{\nu}) \big)\\
 & = 4\big(\psi^3 \phi \big( 3n^2+n(5-4\lambda)+2(\lambda^2-\lambda +1)\big)\\
 & \quad +\lambda \psi^3 x^2 \big( 4n + 2(1-\lambda) \big) \big)
\end{split}
\end{equation}  
so that (\ref{eq5}) vanishes if $(\lambda = \frac{1}{M-1}= \varepsilon_a)$
\begin{equation}\label{eq15}
\gamma + \lambda = 1 \quad \text{and} \quad n=n_{\pm} =(\lambda-\frac{1}{2}) \pm (\lambda+\frac{1}{2})
\end{equation}
while $n_- = -1$ is uninteresting (corresponding to linear shifts of some of the coordinates)
\begin{equation}\label{eq16}
(t^2 - z^2 + \frac{r^2}{M-1})(t+z)^{\frac{2}{M-1}} = C'
\end{equation}
shows that (\ref{eq9}) is not the only possible generalization of the corresponding $M=2$ solutions given in \cite{H94}.
Due to the correspondence found in \cite{BH94} this will give rise to exact solutions
\begin{equation}\label{eq17}
\begin{split}
t-z & = p(\tau, r) = \tau^{a}P(\tau^c r)\\
a + 2&c  + 1   = 0, \quad
P(w)  = \alpha + \beta\frac{w^2}{2}\\
a & = -\frac{M+1}{M-1}, \quad \beta = \frac{-1}{M-1}
\end{split}
\end{equation}
(as well as $a = 2M+1$, $\beta = 1$) of the hydrodynamical equations 
\begin{equation}\label{eq18}
\ddot{p}+ 2(p'\dot{p}' - \dot{p}p'') = \frac{M-1}{r}(2\dot{p}p'+p'^3),
\end{equation}
resp.
\begin{equation*}
\begin{split}
a(a-1)P -\frac{3}{4}(a^2-1)wP' + \frac{(a+1)^2}{4}w^2P'' + (Ma + M-2)P'^2\\
= \frac{M-1}{w}P'^3 + 2a\big( PP'' + (M-1)\frac{PP'}{w}\big).
\end{split}
\end{equation*}

\vspace{1.5cm}
\noindent
\textbf{Acknowledgement.} I am grateful to V.Bach, J.Eggers, and M.Hynek for very valuable discussions.

\end{document}